\documentclass[a4paper,reqno,12pt]{amsart}
\usepackage{amssymb}
\usepackage{ifthen}
\usepackage{graphicx}
\usepackage{color,xcolor}
\usepackage{mathrsfs}
\usepackage{booktabs}
\newtheorem{thm}{Theorem}
\newtheorem{lem}{Lemma}
\newtheorem{cor}{Corollary}
\newtheorem{prop}{Proposition}
\newtheorem{defn}{Definition}
\newtheorem{example}{Example}

\newtheorem{claim}{Claim}

\newtheorem{conj}{Conjecture}
\newtheorem{prob}{Problem}
\newtheorem{rem}{Remark}

\newenvironment{pf}[1][]{%
 \vskip 1mm
 \noindent
 \ifthenelse{\equal{#1}{}}%
  {{\slshape Proof. }}%
  {{\slshape #1.} }%
 }%
{\qed\medskip}

\newcounter{alphabet}
\newcounter{tmp}

\makeatletter
\newcommand{\Ref}[1]{\@ifundefined{r@#1}{}{\setcounter{tmp}{\ref{#1}}\Alph{tmp}}}
\makeatother

\DeclareMathOperator{\RE}{Re}

\newcommand{\ID}{{\mathbb D}}

\def\be{\begin{equation}}
\def\ee{\end{equation}}

\newcommand{\bee}{\begin{enumerate}}
\newcommand{\eee}{\end{enumerate}}

\newcommand{\blem}{\begin{lem}}
\newcommand{\elem}{\end{lem}}
\newcommand{\bthm}{\begin{thm}}
\newcommand{\ethm}{\end{thm}}
\newcommand{\bcor}{\begin{cor}}
\newcommand{\ecor}{\end{cor}}
\newcommand{\beg}{\begin{example}}
\newcommand{\eeg}{\end{example}}
\newcommand{\begs}{\begin{examples}}
\newcommand{\eegs}{\end{examples}}
\newcommand{\bdefe}{\begin{defn}}
\newcommand{\edefe}{\end{defn}}
\newcommand{\bprob}{\begin{prob}}
\newcommand{\eprob}{\end{prob}}
\newcommand{\bques}{\begin{ques}}
\newcommand{\eques}{\end{ques}}
\newcommand{\bei}{\begin{itemize}}
\newcommand{\eei}{\end{itemize}}

\newcommand{\bde}{\begin{deter}}
\newcommand{\ede}{\end{deter}}
\newcommand{\bca}{\begin{case}}
\newcommand{\eca}{\end{case}}
\newcommand{\bcl}{\begin{claim}}
\newcommand{\ecl}{\end{claim}}
\newcommand{\bcon}{\begin{conj}}
\newcommand{\econ}{\end{conj}}
\newcommand{\bcons}{\begin{conjs}}
\newcommand{\econs}{\end{conjs}}
\newcommand{\bprop}{\begin{propo}}
\newcommand{\eprop}{\end{propo}}
\newcommand{\br}{\begin{rem}}
\newcommand{\er}{\end{rem}}
\newcommand{\brs}{\begin{rems}}
\newcommand{\ers}{\end{rems}}
\newcommand{\bo}{\begin{obser}}
\newcommand{\eo}{\end{obser}}
\newcommand{\bos}{\begin{obsers}}
\newcommand{\eos}{\end{obsers}}
\newcommand{\bpf}{\begin{pf}}
\newcommand{\epf}{\end{pf}}
\newcommand{\ba}{\begin{array}}
\newcommand{\ea}{\end{array}}
\newcommand{\beq}{\begin{eqnarray}}
\newcommand{\beqq}{\begin{eqnarray*}}
\newcommand{\eeq}{\end{eqnarray}}
\newcommand{\eeqq}{\end{eqnarray*}}

\newcommand{\ds}{\displaystyle}

\begin{document}
\title[On Harmonic $\nu$-Bloch and $\nu$-Bloch-type mappings]
{On Harmonic $\nu$-Bloch and $\nu$-Bloch-type mappings}


\author[G. Liu]{Gang Liu}
\address{G. Liu, College of Mathematics and Statistics
 (Hunan Provincial Key Laboratory of Intelligent Information Processing and Application),
Hengyang Normal University, Hengyang,  Hunan 421008, People's Republic of China.}
\email{liugangmath@sina.cn}

\author[S. Ponnusamy]{Saminathan Ponnusamy
}
\address{S. Ponnusamy, Stat-Math Unit,
Indian Statistical Institute (ISI), Chennai Centre,
110, Nelson Manickam Road,
Aminjikarai, Chennai, 600 029, India.}
\email{samy@isichennai.res.in, samy@iitm.ac.in}



\begin{abstract}
The aim of this paper is twofold.
One is to introduce the class of harmonic $\nu$-Bloch-type mappings as a generalization of harmonic $\nu$-Bloch mappings
and thereby we generalize some recent results of harmonic $1$-Bloch-type mappings investigated recently by Efraimidis et al. \cite{EGHV}.
The other is to investigate some subordination principles for harmonic Bloch mappings and then
establish Bohr's theorem for these mappings and in a general setting, in some cases.
\end{abstract}

\maketitle \pagestyle{myheadings}
\markboth{G. Liu and S. Ponnusamy}{On Harmonic $\nu$-Bloch and $\nu$-Bloch-type mappings}

\section{Introduction}  \label{sec1}
A significant part of function theory deals with univalent functions, function spaces such as Bloch spaces, Bohr's phenomenon
and their various generalizations. Several authors have contributed a lot to this development, and most importantly,
in the area of planar harmonic mappings. For basic results about harmonic mappings, the reader may refer \cite{CS},
the monograph of Duren \cite{dur2004} and the recent survey of some basic materials from \cite{PR}.
Concerning classical Bloch spaces, see \cite{ACP,BMY,col}. In recent years, Bohr's phenomenon, its various generalizations including
higher dimensional analogues and its harmonic analogues have been studied by various authors.
For more details of the importance, background, development and
results, we refer to the recent survey on this topic \cite{AAP}  and  the references therein.
The recent results on this topic for harmonic mappings may be obtained from \cite{KP,KPS}. Our primary goal here is to
continue to study harmonic Bloch-type mappings and as applications, we consider Bohr's inequality in a general setting.

Throughout we consider complex-valued harmonic mappings in the open unit disk $\mathbb{D}=\{z:\,|z|<1\}$: $\Delta f=f_{z\overline{z}}=0$.
It is well known that every harmonic mapping in $\mathbb{D}$ 
has a {\it canonical decomposition} $f=h+\overline{g}$,
where $h$ and $g$ are analytic functions with $g(0)=0$. Thus, we may express $h$ and $g$ as
\be\label{equ1a}
h(z)=\sum_{n=0}^\infty a_nz^n \quad \text{and}\quad g(z)=\sum_{n=1}^\infty b_nz^n.
\ee
Moreover, the function $f=h+\overline{g}$ is locally univalent and {\it sense-preserving} in $\mathbb{D}$
if and only if its Jacobian $J_f= |f_z|^2-|f_{\overline{z}}|^2=|h'|^2-|g'|^2>0$  in $\mathbb{D}$ by Lewy's theorem (see \cite{lew}),
i.e., $|h'|>|g'|$ or $|\omega_f|<1$ in $\mathbb{D}$, where $\omega_f=g'/h'$ is the dilatation of $f$.

For a given $\nu\in(0,\infty)$, a harmonic mapping $f=h+\overline{g}$ in $\mathbb{D}$
is called a {\it harmonic $\nu$-Bloch mapping} if
\begin{equation*} 
\beta_\nu(f):=\sup_{z\in\mathbb{D}}(1-|z|^2)^\nu(|h'(z)|+|g'(z)|)<\infty.
   \end{equation*}
This defines a seminorm, and the space equipped with the norm
$$||f||_{\mathcal{B}_{H}(\nu)}:=|f(0)|+\beta_\nu(f)
$$
is called the {\it harmonic $\nu$-Bloch space}, denoted by $\mathcal{B}_{H}(\nu)$. It is a Banach space.
In particular the space $\mathcal{B}(\nu)$ defined by
$$\mathcal{B}(\nu)=\{f=h+\overline{g}\in\mathcal{B}_{H}(\nu):\,g\equiv0\}
$$
forms a Banach space equipped with the norm $||f||_{\mathcal{B}(\nu)}:=|f(0)|+\beta_\nu(f)$.
Clearly, $f=h+\overline{g}\in\mathcal{B}_{H}(\nu)$ if and only if $h,g\in\mathcal{B}(\nu)$,
since $$\max\{\beta_\nu(h),\beta_\nu(g)\}\leq\beta_\nu(f)\leq\beta_\nu(h)+\beta_\nu(g).
$$
The harmonic $\nu$-Bloch space $\mathcal{B}_{H}(\nu)$ was introduced in \cite{CPW2011},  which was a generalization of
$\mathcal{B}_{H}(1)$ that was studied by Colonna  in \cite{col} as a generalization of classical Bloch space $\mathcal{B}(1)$.
One can refer to \cite{ACP,BMY,CGH,CPW2012,pom,zjf} for information on $\mathcal{B}(1)$ and its extension.
Motivated by results on analytic Bloch functions, Efraimidis et at. \cite{EGHV} introduced  harmonic Bloch-type mappings,
which coincide with the following harmonic $1$-Bloch-type mappings.

\bdefe
For a given $\nu\in(0,\infty)$, a harmonic mapping $f$ on $\mathbb{D}$
is called a {\it harmonic $\nu$-Bloch-type mapping} if
\begin{equation*}
\beta^*_\nu(f):=\sup_{z\in\mathbb{D}}(1-|z|^2)^\nu\sqrt{|J_f(z)|}<\infty.
\end{equation*}
We write $\mathcal{B}^*_{H}(\nu)$ for the space of all such mappings and we call
$$||f||_{\mathcal{B}^*_{H}(\nu)}:=|f(0)|+\beta^*_\nu(f),
$$
the pseudo-norm of $f$.
\edefe

In Section \ref{sec2}, we will see that $\mathcal{B}^*_{H}(\nu)$ is not a linear space for any $\nu>0$.
Because
$$(1-|z|^2)^\nu\sqrt{|J_f(z)|}\leq(1-|z|^2)^\nu(|f_z(z)|+|f_{\overline{z}}(z)|),\quad z\in\mathbb{D},
$$
it is clearly that $\mathcal{B}_{H}(\nu)\subset \mathcal{B}^*_{H}(\nu)$ and thus, the space $\mathcal{B}^*_{H}(\nu)$ is a
generalization of $\mathcal{B}_{H}(\nu)$.
In addition, in the case of analytic functions $f$, these spaces coincide and thus, we have $$||f||_{\mathcal{B}^*_{H}(\nu)}=||f||_{\mathcal{B}_{H}(\nu)}=||f||_{\mathcal{B}(\nu)}.
$$

One of the aims of this article is to generalize some of the known results of harmonic $\nu$-Bloch mappings and
$\nu$-Bloch-type mappings (especially, the results of \cite{EGHV}). The paper is divided into sections as follow:
In Section \ref{sec2}, for the function spaces $\mathcal{B}_{H}(\nu)$ and $\mathcal{B}^*_{H}(\mu)$,
we investigate  its  {\it affine and linear invariance}, and the inclusion relations under particular conditions.
In Section \ref{sec3}, we find a connection between these function spaces
and the space of {\it uniformly locally univalent} harmonic mappings.
Moreover, some {\it subordination} principles concerning the spaces $\mathcal{B}_{H}(1)$ and $\mathcal{B}^*_{H}(1)$ are also investigated.
In Section \ref{sec4}, we give the growth and coefficients estimates for sense-preserving mappings in $\mathcal{B}^*_{H}(\nu)$.
Finally, as applications of our investigations, we determine the {\it Bohr radius} for functions in $\mathcal{B}(\nu)$,
and {\it $p$-Bohr radius} for  functions in $\mathcal{B}_{H}(\nu)$ and $\mathcal{B}^*_{H}(\nu)$ (sense-preserving)
in Section \ref{sec5}.

\section{Affine and linear invariance and inclusion relations} \label{sec2}

Throughout the article $\nu$ is a constant in the interval $(0,\infty)$.
We first discuss the affine and linear invariance of $\mathcal{B}_{H}(\nu)$ and $\mathcal{B}^*_{H}(\nu)$.
Let $L$ be a family of harmonic mappings defined in $\mathbb{D}$. Then the family
$L$ is said to be affine invariant if $A\circ f\in L$ for each $f\in L$ and for all affine mappings $A$ of the form
$A(z)=az+b\overline{z}$ ($a,b\in\mathbb{C}$). The family $L$ is called linear invariant if  for each $f\in L$,
$$f\circ \varphi_\alpha\in L \quad \forall ~~\varphi_\alpha(z)
=\frac{z+\alpha}{1+\overline{\alpha}z}\in\rm{Aut}(\mathbb{D}).
$$

\begin{prop}    \label{prop1}
\begin{enumerate}
\item Both $\mathcal{B}_{H}(\nu)$ and $\mathcal{B}^*_{H}(\nu)$ are affine invariant.
\item  Each of $\mathcal{B}(\nu)$,  $\mathcal{B}_{H}(\nu)$ and $\mathcal{B}^*_{H}(\nu)$ is linear invariant.
\end{enumerate}
\end{prop}

\bpf
(1) Let $f=h+\overline{g}$ and $A(z)=az+b\overline{z}$ ($a,b\in\mathbb{C}$).
Then
$$A\circ f=ah+bg+\overline{\overline{a}g+\overline{b}h}
$$
and thus,
$$|(ah+bg)'|+|(\overline{a}g+\overline{b}h)'|\leq(|a|+|b|)(|h'|+|g'|)
$$
and $$J_{A\circ f}=|(ah+bg)'|^2-|(\overline{a}g+\overline{b}h)'|^2=(|a|^2-|b|^2)J_f.$$
The desired conclusion now easily follows.

(2) We only need to prove that $\mathcal{B}_H^*(\nu)$ is linear invariant.
For $\varphi_\alpha(z)\in\rm{Aut}(\mathbb{D})$, we have
$$J_{f\circ \varphi_\alpha}(z)=|\varphi'_\alpha(z)|^2J_f(\varphi_\alpha(z))
~\mbox{ and }~ 1-|\varphi_\alpha(z)|^2=|\varphi'_{\alpha}(z)|(1-|z|^2),
$$
and obtain
\begin{align*}
(1-|z|^2)^\nu\sqrt{|J_{f\circ \varphi_\alpha}(z)|}&=(1-|z|^2)^\nu|\varphi'_\alpha(z)|\sqrt{|J_f(\varphi_\alpha(z))|}\\
&=\left(\frac{1-|z|^2}{1-|\varphi_\alpha(z)|^2}\right)^{\nu-1}(1-|\varphi_\alpha(z)|^2)^v\sqrt{|J_f(\varphi_\alpha(z))|}\\
&\leq\left(\frac{|1+\overline{\alpha} z|^2}{1-|\alpha|^2}\right)^{\nu-1}\beta^*_\nu(f)\\
&\leq \left(\frac{1+|\alpha|}{1-|\alpha|}\right)^{|\nu-1|}\beta^*_\nu(f),\quad z\in\mathbb{D}.
\end{align*}
Now it is obvious that if $f\in\mathcal{B}^*_{H}(\nu)$, then $f\circ \varphi_\alpha\in\mathcal{B}^*_{H}(\nu)$
for each $\varphi_\alpha\in\rm{Aut}(\mathbb{D})$.

Similarly, $\mathcal{B}(\nu)$ is linear invariant so that $\mathcal{B}_{H}(\nu)$ is also linear invariant.
\epf

For each $\nu>0$, although both $\mathcal{B}(\nu)$ and $\mathcal{B}_H(\nu)$ are Banach spaces,
the following example shows that $\mathcal{B}^*_{H}(\nu)$ is not a linear space.
It also shows that some functions in $\mathcal{B}^*_{H}(\nu)$ may grow arbitrarily fast.
Therefore, to study certain properties of functions in $\mathcal{B}^*_{H}(\nu)$ in what follows,
we shall restrict harmonic mappings  to be sense-preserving.

\begin{example} \label{exa1}
{\rm
Let $f=h+\overline{h}$, where $h(z)=(\mu-1)^{-1}(1-z)^{1-\mu}$ for some $\mu>2\nu+1$.
Clearly, we have $f(z)$ and the identity function $z$ belong to $\mathcal{B}^*_{H}(\nu)$ whereas $F(z)=f(z)+z$ does not, since
$$(1-x^2)^{2\nu}|J_{F}(x)|=(1+x)^{2\nu}\frac{2+(1-x)^\mu}{(1-x)^{\mu-2\nu}}\rightarrow\infty ~\mbox{ as }~ (0,1)\ni x\rightarrow1^-.
$$
}
\end{example}

Next we deal with the inclusion relations $\mathcal{B}(\nu)\subset\mathcal{B}_{H}(\nu)\subseteq\mathcal{B}^*_{H}(\nu)$.

\begin{prop}    \label{prop2}
Let $\mu$ and $\nu$ be two constants with $0<\mu<\nu$. We have
\begin{enumerate}
\item $\mathcal{B}(\nu)\subset\mathcal{B}_{H}(\nu)\subset\mathcal{B}^*_{H}(\nu)$;
\item $\mathcal{B}(\mu)\subset\mathcal{B}(\nu)$,  $\mathcal{B}_{H}(\mu)\subset\mathcal{B}_{H}(\nu)$ and $\mathcal{B}^*_{H}(\mu)\subset\mathcal{B}^*_{H}(\nu)$.
\end{enumerate}
\end{prop}

\bpf
(1) It only needs to find a function $f_\nu\in\mathcal{B}^*_{H}(\nu)\backslash\mathcal{B}_{H}(\nu)$ for each $\nu>0$.
For the sake of the context later, we will prove that the one parameter family of functions $f_{\nu,t}\in\mathcal{B}^*_{H}(\nu)\backslash\mathcal{B}_{H}(\nu)$
for each $\nu>0$, where
\begin{equation}\label{equ1}
f_{\nu,t}(z)=h_{\nu}(z)+ \overline{g_{\nu,t}(z)}, \quad t\in[0,1),
\end{equation}
with
\begin{equation}\label{equ2}
h_\nu(z)=
\begin{cases}
-\log(1-z) & \mbox{for }~  \nu=1/2,\\
(\nu-1/2)^{-1}\left[(1-z)^{1/2-\nu}-1\right] & \mbox{for }~ \nu\neq1/2,
\end{cases}
\end{equation}
and
\begin{align*}
g_{\nu,t}(z)=
\begin{cases}
-\log(1-z)-(1-t)z  &\mbox{for }~  \nu=1/2,\\
\ds \frac{z}{1-z}+(1-t)\log(1-z) & \mbox{for }~ \nu=3/2,\\
\ds \frac{(1-z)^{1/2-\nu}-1}{\nu-1/2}-(1-t)\frac{(1-z)^{3/2-\nu}-1}{\nu-3/2} & \mbox{for }~   \nu\neq 1/2, 3/2.
\end{cases}
\end{align*}
Fix $t\in[0,1)$. A direct computation reveals that $f_{\nu,t}$ is sense-preserving in $\mathbb{D}$
with the dilatation $\omega_{f_{\nu,t}}(z)=t+(1-t)z$ for each $\nu>0$.
Again, by computation, we have
\begin{align*}
(1-|z|^2)^\nu\sqrt{|J_{f_{\nu,t}}(z)|}&=\frac{(1-|z|^2)^\nu}{|1-z|^{\nu+1/2}}|\sqrt{1-|\omega_{f_{\nu,t}}(z)|^2}\\
&=(1+|z|)^\nu\frac{(1-|z|)^{\nu}}{|1-z|^{\nu}}\sqrt{\frac{1-|z|^2-2t\RE(\overline{z}(1-z))-t^2|1-z|^2}{|1-z|}}\\
&\leq2^{\nu+1/2}\sqrt{1+t},
\end{align*}
which gives $f_{\nu,t}\in\mathcal{B}^*_H(\nu)$ for each $\nu>0$.
Since
$$(1-x^2)^\nu|h'_\nu(x)|=\frac{(1+x)^\nu}{\sqrt{1-x}}\rightarrow\infty~\mbox{ as }~ (0,1)\ni x\rightarrow1^-,
$$
we obtain that for each $\nu>0$, $h_\nu\not\in\mathcal{B}_{H}(\nu)$ and thus, $f_{\nu,t}\not\in\mathcal{B}_{H}(\nu)$.

(2)  Let $0<\mu<\nu$. Clearly, $\mathcal{B}(\mu)\subseteq\mathcal{B}(\nu)$,  $\mathcal{B}_{H}(\mu)\subseteq\mathcal{B}_{H}(\nu)$ and $\mathcal{B}^*_{H}(\mu)\subseteq\mathcal{B}^*_{H}(\nu)$.
Then the inclusions $\mathcal{B}_{H}(\mu)\subset\mathcal{B}_{H}(\nu)$ and $\mathcal{B}^*_{H}(\mu)\subset\mathcal{B}^*_{H}(\nu)$
obviously follow by (1) if we prove $\mathcal{B}(\mu)\subset\mathcal{B}(\nu)$.
For this, we simply consider the function $f_\nu$ satisfying $f'_\nu(z)=(1-z)^{-\nu}$, it is easy to see that $f_\nu\in\mathcal{B}(\nu)\backslash\mathcal{B}(\mu)$.
This completes the proof.
\epf

It is natural to ask for the structure of the set $\mathcal{B}^*_H(\nu)\backslash\mathcal{B}_{H}(\nu)$.

\begin{prop}  \label{prop3}
Let $f=h+\overline{g}$ be a harmonic mapping in $\mathbb{D}$.
Then $f\in\mathcal{B}_{H}(\nu)$ if and only if
$f\in\mathcal{B}^*_{H}(\nu)$ and either $h\in\mathcal{B}(\nu)$ or $g\in\mathcal{B}(\nu)$.
We get
$$\mathcal{B}^*_H(\nu)\backslash\mathcal{B}_{H}(\nu)=\{f=h+\overline{g}\in\mathcal{B}^*_H(\nu):\,
h\not\in\mathcal{B}(\nu)~\text{and}~g\not\in\mathcal{B}(\nu)\}.
$$
\end{prop}
\bpf It suffices to observe that
$|h'|\leq\sqrt{|J_f|}+|g'|$  and $| g'|\leq\sqrt{|J_f|}+|h'|$
for a harmonic mapping $f=h+\overline{g}$.
\epf

The following question arises.

\bprob
Suppose that $f\in\mathcal{B}^*_{H}(\nu)$. Does there exist a constant $c(\nu)$ depending only on $\nu$
such that $f\in\mathcal{B}_{H}(c(\nu))$?
\eprob

In order to give an affirmative answer to this problem,  we need  some extra conditions based on the following observation
for the function
$f=h+\overline{h}$, where $h(z)=\exp \left ( (1+z)/(1-z)\right ).$
Clearly, $f\in\mathcal{B}^*_{H}(\nu)$ for all $\nu>0$.
However,  $f\not\in\mathcal{B}_{H}(\nu)$ for any $\nu>0$, since  $h\not\in\mathcal{B}(\nu)$,
which can be deduced from
$$(1-x^2)^{\nu}|h'(x)|=2(1+x)^{2\nu-2}
\left[\left(\frac{1-x}{1+x}\right)^{\nu-2}e^{\frac{1+x}{1-x}}\right]\rightarrow\infty~\mbox{ as }~(0,1)\ni x\rightarrow1^-.
$$


\begin{prop}   \label{prop4}
Let $f$ be a locally univalent  harmonic mapping in $\mathbb{D}$.
If $f\in\mathcal{B}^*_{H}(\nu)$, then $f\in\mathcal{B}_{H}(\nu+1/2)$.
Moreover, the constant $1/2$ is sharp for each $\nu>0$.
\end{prop}
\bpf
Note that $f\in\mathcal{B}_H(\nu)$ (resp. $\mathcal{B}^*_H(\nu)$) if and only if $\overline{f}\in\mathcal{B}_H(\nu)$ (resp. $\mathcal{B}^*_H(\nu)$).
Without loss of generality, we may thus assume that $f=h+\overline{g}$  is sense-preserving with the dilatation $\omega=\omega_f$ so that
$$ g'=\omega h' ~\mbox{and }~ J_f= |h'|^2 (1-|\omega|^2) ~\mbox{ or }~ |h'|=\sqrt{\frac{J_f}{1-|\omega|^2}}.
$$
It follows (see \cite[Corollary~1.3]{gar}) that
$$|\omega(z)|\leq \frac{|z|+|\omega(0)|}{1+|\omega(0)z|}, \quad  z\in \mathbb{D}.
$$
Now we suppose that $f\in\mathcal{B}^*_{H}(\nu)$. Then we get
$$(1-|z|^2)^\nu\sqrt{|J_f(z)|}\leq\beta^*_\nu(f)<\infty , \quad  z\in \mathbb{D}.
$$
Consequently, we have $|g'(z)|< |h'(z)|$ in $\ID$, where
\begin{align}
 |h'(z)|=&\sqrt{\frac{J_f(z)}{1-|\omega(z)|^2}}\nonumber\\
\leq &\frac{\beta^*_\nu(f)}{(1-|z|^2)^\nu} \frac{1}{\sqrt{1-|\omega(z)|^2}} \label{equ3}\\
\nonumber \leq &\frac{\beta^*_\nu(f)}{(1-|z|^2)^\nu}\left[1-\left(\frac{|z|+|\omega(0)|}{1+|\omega(0)||z|}\right)^2\right]^{-1/2}\\
\nonumber = &\frac{\beta^*_\nu(f)}{(1-|z|^2)^\nu}\left[\frac{1+|\omega(0)|\,|z|}{\sqrt{(1-|z|^2)(1-|\omega(0)|^2)}}\right]\\
\leq &\frac{\beta^*_\nu(f)}{(1-|z|^2)^{\nu+\frac{1}{2}}}\sqrt{\frac{1+|\omega(0)|}{1-|\omega(0)|}} \label{equ4},
\end{align}
which shows that $h$ (and hence $g$) belongs to $\mathcal{B}_{H}(\nu+1/2)$. Hence, $f\in\mathcal{B}_{H}(\nu+1/2)$.

To see that the constant $1/2$ is sharp for each $\nu>0$, it suffices to check for the function $f_{\nu,0}=h_\nu+\overline{g_{\nu,0}}$ defined by \eqref{equ1}.
From the proof of Proposition \ref{prop2},  the function $f_{\nu,0}~(\in\mathcal{B}^*_{H}(\nu))$ is sense-preserving in $\mathbb{D}$.
On the other hand, it is easy to see that $h_\nu\in\mathcal{B}(\nu+1/2)$, which implies $g_{\nu,0}\in\mathcal{B}(\nu+1/2)$ and thus, $f_\nu\in\mathcal{B}_H(\nu+1/2)$. However, we have that for any $0<\varepsilon<\nu+1/2$, $h_\nu\not\in\mathcal{B}(\varepsilon)$,
which means  $f_{\nu,0}\not\in\mathcal{B}_H(\varepsilon)$.
We complete the proof.
\epf

\section{Uniformly locally univalent and subordination principles} \label{sec3}

\subsection{Connection with uniformly locally univalent harmonic mappings} \label{sec3.1}

Motivated by the characterization of Bloch space $\mathcal{B}(1)$ and the recent work of the authors \cite{LS} concerning
equivalent conditions of {\it uniformly locally univalent} (briefly, ULU) harmonic mappings,
we will show the connections among harmonic $\nu$-Bloch, $\nu$-Bloch-type  mappings and ULU harmonic mappings.

We first introduce the notion and some properties of  ULU harmonic mappings.
A harmonic  mapping $f=h+\overline{g}$ in $\mathbb{D}$ is called ULU if there exists a constant
$\rho>0$ such that $f$ is univalent on the hyperbolic disk
$$D_h(a,\rho)=\left \{z\in\mathbb{D}:\,\left |\frac{z-a}{1-\overline{a}z}\right |<\tanh\rho \right \},
$$
of radius $\rho$, for every $a\in\mathbb{D}$.
One of equivalent conditions of ULU is stated in terms of the {\it pre-Schwarzian derivative or norm}.
Let $f$ be a locally univalent harmonic mapping in $\mathbb{D}$.
The pre-Schwarzian derivative and the norm of $f$ are defined as \cite{HM2015} (see also \cite{CDO})
\begin{equation*}
P_f=(\log J_f)_z,\quad z\in\mathbb{D},
\quad \text{and}\quad  ||P_{f}||=\sup_{z\in\mathbb{D}}(1-|z|^2)|P_{f}(z)|,
\end{equation*}
respectively.
Clearly, the two definitions coincide with the corresponding definitions in the analytic case.
Similar to the proof of \cite[Theorem~7]{HM2015}, the function $f=h+\overline{g}$ in $\mathbb{D}$ is ULU if and only if $||P_f||<\infty$
(see also \cite[Theorem 4.1]{LS}). Several equivalent conditions of ULU mappings can be found in these three papers and the references therein.

Now let's restrict $f$ to be analytic in $\mathbb{D}$. It is well-known that $f\in\mathcal{B}(1)$ if and only if
there exists a constant $c>0$ and a univalent analytic function $F$ such that $f=c\log F'$ (see \cite{pom}).
On the other hand,  $f$ is ULU if and only if there exists a constant $c>0$ and a univalent analytic
function $F$ such that $f'=(F')^c$ (see \cite[Theorem~2]{yam1977}). Thus, $f\in\mathcal{B}(1)$ if and only if
there exists a ULU analytic function $F$ such that $f=\log F'$.
Furthermore, a harmonic mapping $f=h+\overline{g}$ belongs to $\mathcal{B}_H(1)$
if and only if there exist two ULU analytic functions $H$ and $G$ such that $f=\log H'+\overline{\log G'}$.
A natural question is to ask: {\it What about the characterization of $\mathcal{B}^*_H(1)$?}
The following theorem and example show some extraneous complexities of the structure of
the space $\mathcal{B}^*_H(1)$, which are different from Example \ref{exa1}.

\bthm \label{thm1}
Let $F=H+\overline{G}$ be sense-preserving and {\rm ULU} in $\mathbb{D}$. Then for each $\varepsilon\in\overline{\mathbb{D}}$, the function $f_\varepsilon=h_\varepsilon+\overline{g_\varepsilon}$ belongs to $\mathcal{B}^*_H(1)$, where $h_\varepsilon=\log(H'+\varepsilon G')$
and $\omega=g'_\varepsilon/h'_\varepsilon$ is bounded in $\ID$.
\ethm
\bpf
Suppose that $F=H+\overline{G}$ is a sense-preserving and ULU in $\mathbb{D}$. It follows
from \cite[Theorem~4.1]{LS} that $||P_{H+\varepsilon G}||<\infty$
for all $\varepsilon\in\overline{\mathbb{D}}$.
By assumption, for each $\varepsilon\in\overline{\mathbb{D}}$, we have
\begin{align*}
(1-|z|^2)\sqrt{|J_{f_\varepsilon}(z)|}&
\leq(1-|z|^2)|h'_\varepsilon(z)|(1+\sup_{z\in\mathbb{D}}|\omega(z)|)\\
&=(1-|z|^2)\left|\frac{H''+\varepsilon G''}{H'+\varepsilon G'}\right|(1+\sup_{z\in\mathbb{D}}|\omega(z)|)\\
&\leq||P_{H+\varepsilon G}||(1+\sup_{z\in\mathbb{D}}|\omega(z)|)<\infty
\end{align*}
and the assertion follows.
\epf

\begin{example} \label{exa2}
{\rm
Consider the function $f=h+\overline{g}$ in $\mathbb{D}$ with the dilatation $\omega_f(z)=e^{i\theta}z$,
where
$$h=\log H' ~\mbox{ and }~ H(z)=\exp \left ({\sqrt{\frac{1+z}{1-z}}} \right )=:\exp (q(z)), \quad z\in\ID,
$$
and the principal branch of the square root is chosen such that $q(0)=1$.
We claim that $f\in\mathcal{B}^*_H(1)\backslash\mathcal{B}_H(1)$ and $H$ is locally univalent but not \textrm{ULU} in $\mathbb{D}$.
To do this, straightforward  computations give that
$$H'(z)=\frac{1}{(1-z)^2}\sqrt{\frac{1-z}{1+z}}\exp \left ({\sqrt{\frac{1+z}{1-z}}} \right )\neq0,\quad z\in\mathbb{D},
$$
and
$$||P_{H}||=\sup_{z\in\mathbb{D}}(1-|z|^2)\left|\frac{(1+2z)\sqrt{1-z}+\sqrt{1+z}}{(1-z^2)\sqrt{1-z}}\right|=\infty,
$$
showing that $H$ is locally univalent but not ULU in $\mathbb{D}$. 
Again, elementary computations show that
\begin{align*}
(1-|z|^2)\sqrt{|J_{f}(z)|}&
=(1-|z|^2)|h'(z)|\sqrt{1-|z|^2}\\
&=(1-|z|^2)\left|\frac{(1+2z)\sqrt{1-z}+\sqrt{1+z}}{(1-z^2)\sqrt{1-z}}\right|\sqrt{1-|z|^2}\\
&\leq\left|(1+2z)\sqrt{1-z}+\sqrt{1+z}\right|\sqrt{1+|z|}<\infty,\quad z\in\mathbb{D},
\end{align*}
which implies  $f\in\mathcal{B}^*_H(1)$.
Moreover,   because  $H$ is not ULU, we find that $h\not\in\mathcal{B}(1)$ and thus, $f\not\in\mathcal{B}_H(1)$.
Hence we conclude that, $f\in\mathcal{B}^*_H(1)\backslash\mathcal{B}_H(1)$.
}
\end{example}

Although Theorem \ref{thm1} is a generalization of \cite[Theorem~2]{EGHV},
we can't give a complement characterization of $\mathcal{B}^*_H(1)$, let alone to $\mathcal{B}^*_H(\nu)$ ($\nu>0$).
However, we will see that any ULU harmonic mapping is a $\nu$-Bloch-type mapping for some $\nu>0$.
Recall that $f$ is ULU if and only if $||P_f||<\infty$.
In view of this, to describe our result more precisely, we define the set
$$\mathbb{B}_H(\nu)=\{f:~f ~\text{is ~a~ locally~ univalent~ harmonic~mapping~in}~\mathbb{D}~\text{with}~||P_f||\leq\nu\}
$$
and its subset  $\mathbb{B}(\nu)$ of all analytic functions in $\mathbb{B}_H(\nu)$.

\bthm  \label{thm2}
For any $\nu>0$, we have  $\mathbb{B}_H(\nu)\subset\mathcal{B}^*_H(\nu/2)$.
In particular, $\mathbb{B}(\nu)\subset\mathcal{B}(\nu/2)$.
Moreover, these two inclusions are best possible.
\ethm
\bpf
Assume $f\in\mathbb{B}_H(\nu)$ for some $\nu>0$. Note that $||P_f||=||P_{\overline{f}}||$.
Without loss of generality, we may assume that $f$ is sense-preserving. Then, because $f_z(0)\neq 0$, we may consider
$$F(z)=\frac{f(z)-f(0)}{f_z(0)}.
$$
Then $F$ is sense-preserving in $\mathbb{D}$ with the normalization $F(0)=F_z(0)-1=0$. We have $||P_F||=||P_f||$ and thus, $F\in\mathbb{B}_H(\nu)$.
It follows from \cite[Theorem~6.1]{LS} that
$$J_F(z)\leq(1-|F_{\overline{z}}(0)|^2)\left(\frac{1+|z|}{1-|z|}\right)^\nu,\quad z\in\mathbb{D},
$$
which implies $F\in\mathcal{B}^*_H(\nu/2)$.
Since $\mathcal{B}^*_H(\nu)$ preserves affine invariance for each $\nu>0$, we get $f\in\mathcal{B}^*_H(\nu/2)$.
Clearly, $\mathbb{B}_H(\nu)\subset\mathcal{B}^*_H(\nu/2)$ from Example \ref{exa1}.
The sharpness follows if we choose
\begin{eqnarray*}
f(z)=f_{\nu}(z)=\int_0^z\left(\frac{1+t}{1-t}\right)^{\nu/2}dt
+\overline{b_1\int_0^z\left(\frac{1+t}{1-t}\right)^{\nu/2}dt},\quad  z\in\mathbb{D},
\end{eqnarray*}
where $|b_1|<1$. Indeed, it is easy to see that $||P_{f_\nu}||=\nu$ and $f_\nu\in\mathcal{B}^*_H(\nu/2)$ but
$f_\nu\not\in\mathcal{B}^*_H(\varepsilon)$ for any $0<\varepsilon<\nu/2$.

If $f$ is restricted to be analytic, then a similar proof shows that $\mathbb{B}(\nu)\subset\mathcal{B}(\nu/2)$.
The sharpness can be easily seen by considering the above function $f_\nu$ with $b_1=0$.
\epf


\subsection{Subordination principles} \label{sec3.2}
Every bounded analytic function in $\mathbb{D}$ belongs to the (analytic) Bloch space $B(1)$.
This fact also holds for harmonic mappings (see  \cite[Theorem~3]{col}) and for a simpler proof
of it (using subordination), we refer to  \cite[Theorem~A]{CPW2011}. That is, if a harmonic mapping
$f$ is bounded in $\mathbb{D}$, then $f$ belongs to the (harmonic) Bloch space $\mathcal{B}_H(1)$.
Next we will investigate some subordination principles for some harmonic Bloch mappings.


Let $\mathcal{A}_D$ denotes the class of analytic functions $\phi:\mathbb{D}\rightarrow\mathbb{D}$
and  $\mathcal{A}_D^0$ denotes the subclass of $\mathcal{A}_D$ with the normalization  $\phi(0)=0$.
In 2000, Schaubroeck  \cite{sch} generalized the notion of subordination from analytic functions to harmonic mappings.
Let $f$ and $F$ be two harmonic mappings in  $\mathbb{D}$.
Then $f$ is {\it subordinate}  to $F$,  denoted
by $f\prec F$, if there is a function $\phi\in\mathcal{A}_D^0$  such that $f=F\circ \phi$.
We denote $f\preceq F$ if there exists a function $\phi\in\mathcal{A}_D$  such that $f=F\circ \phi$.
Clearly, if $f\prec F$ then $f\preceq F$.


\bthm\label{thm3} {\rm (Subordination principle)}
Let $f$ and $F$ be two harmonic mappings in $\mathbb{D}$.
If $f\preceq F$ and  $F\in\mathcal{B}_H(1)$ (resp. $\mathcal{B}^*_H(1)$), then $f\in\mathcal{B}_H(1)$ (resp. $\mathcal{B}^*_H(1)$).
In particular, if $f\preceq F$ and $F\in\mathcal{B}(1)$, then $f\in\mathcal{B}(1)$.
\ethm
\bpf
We just need to prove the case of $\nu$-Bloch-type mappings since the proof of the remaining cases are similar.
Assume that $f\preceq F$ and $F\in\mathcal{B}^*_H(1)$. Then there exists a function $\phi\in\mathcal{A}_D$ such
that $f=F\circ \phi$. We find that
$$J_{f}(z)=J_F(\phi(z)) |\phi'(z)|^2
$$
and by the Schwarz-Pick lemma, we get $(1-|z|^2)|\phi'(z)|\leq1-|\phi(z)|^2$. Consequently,
\begin{align*}
(1-|z|^2)\sqrt{|J_{f}(z)|}& =(1-|z|^2)|\phi'(z)|\sqrt{|J_F(\phi(z))|}\\
&\leq(1-|\phi(z)|^2)\sqrt{|J_F(\phi(z))|}\leq\beta^*_1(F)<\infty,\quad z\in\mathbb{D},
\end{align*}
which clearly shows that $f\in\mathcal{B}^*_H(1)$.
\epf

\br \label{rem1}
{\rm We remind that  $f=h+\overline{g}\in\mathcal{B}_H(1)$ does not mean that either $h$, $g$ or $f$ is bounded
even if $f$ is sense-preserving in $\mathbb{D}$. For instance, consider
$$f_1(z)=h(z)+\overline{g(z)}=\log(1-z)+\overline{z+\log(1-z)}=\overline{z}+2\log|1-z|.
$$
and
$$f_2(z)=h(z)-\overline{g(z)}=\log(1-z)-\overline{z+\log(1-z)}=-\overline{z}+2i\arg(1-z).
$$
Then it is easy to verify that $f_1,f_2\in\mathcal{B}_H(1)$, and both $f_1$ and $f_2$ are sense-preserving in $\mathbb{D}$.
However,   except $f_2$, neither $h$, nor $g$ nor $f_1$ is bounded in $\mathbb{D}$.
}\er

 \section{Growth and coefficients estimates} \label{sec4}

In this section, we investigate some growth and coefficients estimates for functions in $\mathcal{B}^*_H(\nu)$.
For corresponding results in the case of $\mathcal{B}_H(\nu)$, the reader can refer to \cite{CPW2011,zjf}.

\bthm \label{thm4}
Suppose that $f=h+\overline{g}\in\mathcal{B}^*_H(\nu)$ is sense-preserving in $\mathbb{D}$ with the dilatation $\omega_f$,
where $h$ and $g$ are given by \eqref{equ1a}. Then
$$\max\{|h(z)-a_0|,|g(z)|\}\leq\beta^*_\nu(f)\sqrt{\frac{1+|\omega_f(0)|}{1-|\omega_f(0)|}}h_\nu(r),\quad |z|=r<1,$$
where $h_\nu$ is defined by \eqref{equ2}.
The estimate is sharp in order of magnitude for each  $\nu>1/2$.
If $\nu<1/2$, then each of $h,g,f$ is bounded in $\mathbb{D}$.
\ethm
\bpf
Let $|z|=r<1$. 
Following the proof of  Proposition \ref{prop4} and  \eqref{equ4}, because $f$ is sense-preserving, we have
\begin{align*}
\max\{|h(z)-a_0|,|g(z)|\}&=\max\left\{\left|z\int_0^1h'(tz)dt\right|,\left|z\int_0^1g'(tz)dt\right|\right\}\\
&\leq r\int_0^1|h'(tz)|dt\\
&\leq\beta^*_\nu(f)\sqrt{\frac{1+|\omega_f(0)|}{1-|\omega_f(0)|}}\int_0^1\frac{r}{(1-r^2t^2)^{\nu+1/2}}dt\\
&\leq\beta^*_\nu(f)\sqrt{\frac{1+|\omega_f(0)|}{1-|\omega_f(0)|}}\int_0^1\frac{r}{(1-rt)^{\nu+1/2}}dt\\
&=\beta^*_\nu(f)\sqrt{\frac{1+|\omega_f(0)|}{1-|\omega_f(0)|}}h_\nu(r).
\end{align*}

For each $\nu>1/2$, the sharpness of the order of magnitude can be seen from the functions
$f_{\nu,t}=h_\nu+\overline{g_{\nu,t}}$ defined by \eqref{equ1} for $t\in[0,1)$.
Clearly, it is sharp for $h_\nu$ from its formulation.
Fix $t\in[0,1)$. It is also sharp for the function $g_{\nu,t}$, since for $x\in(0,1)$ and any $\varepsilon>0$,
\begin{align*}
&(1-x^2)^{\nu-1/2-\varepsilon}|g_{\nu,t}(x)|\\
=&\frac{(1+x)^{\nu-1/2-\varepsilon}}{(1-x)^{\varepsilon}}
\left|\frac{1-(1-x)^{\nu-1/2}}{\nu-1/2}-\frac{1-t}{\nu-3/2}\left[(1-x)-(1-x)^{\nu-1/2}\right]\right|\rightarrow\infty
\end{align*}
as $x\rightarrow1^{-}$ when $\nu>1/2$ but $\nu\neq3/2$, and
$$(1-x^2)^{1-\varepsilon}|g_{3/2,t}(x)|=(1+x)^{1-\varepsilon}/(1-x)^{\varepsilon}|x+(1-t)(1-x)\log(1-x)|\rightarrow\infty
$$
 as $x\rightarrow1^{-}$.

If $\nu<1/2$, then for $|z|=r<1$ we have
\begin{align*}
\max\{|h(z)-a_0|,|g(z)|\}\leq&\beta^*_\nu(f)\sqrt{\frac{1+|\omega_f(0)|}{1-|\omega_f(0)|}}(1/2-\nu)^{-1}(1-(1-r)^{1/2-\nu})\\
\leq&\beta^*_\nu(f)\sqrt{\frac{1+|\omega_f(0)|}{1-|\omega_f(0)|}}(1/2-\nu)^{-1}.
\end{align*}
Obviously, both $h$ and $g$ are bounded in $\mathbb{D}$ and thus, $f$ is also bounded in $\mathbb{D}$.
\epf

If $f\in\mathcal{B}^*_H(1/2)$ is sense-preserving in $\mathbb{D}$,
then the boundedness of $f$ is uncertain, which may be verified easily by considering
the two functions  $f_1$ and $f_2$ in Remark \ref{rem1}.
Indeed,  $f_1,~f_2\in\mathcal{B}^*_H(1/2)$.

\bthm \label{thm5}
Suppose that $f=h+\overline{g}\in\mathcal{B}^*_H(\nu)$ is sense-preserving in $\mathbb{D}$ with the dilatation $\omega_f$,
where $h$ and $g$ are given by \eqref{equ1a}. Then
$$|b_1|<|a_1|\leq\frac{\beta^*_\nu(f)}{\sqrt{1-|\omega_f(0)|^2}},
$$
and
$$\max\{|a_n|,|b_n|\}\leq\beta^*_\nu(f)\left(\frac{e}{2\nu+1}\right)^{\nu+1/2}
\sqrt{\frac{1+|\omega_f(0)|}{1-|\omega_f(0)|}}(n+2\nu)^{\nu-1/2},\quad n\geq2.
$$
\ethm

\bpf
The first inequality follows if we set $z=0$ in \eqref{equ3}.
For the second inequality, we recall from \eqref{equ4} that
\beq  \nonumber
|g'(z)|^2<|h'(z)|^2\leq \frac{1+|\omega_f(0)|}{1-|\omega_f(0)|}\frac{\beta^*_\nu(f)^2}{(1-|z|^2)^{2\nu+1}},\quad z\in\mathbb{D}.
\eeq
We integrate this inequality over the circle $|z|=r$ and get
\begin{equation*}
\sum_{n=1}^\infty n^2|b_n|^2r^{2(n-1)}<\sum_{n=1}^\infty n^2|a_n|^2r^{2(n-1)}
\leq \frac{1+|\omega_f(0)|}{1-|\omega_f(0)|}\frac{\beta^*_\nu(f)^2}{(1-r^2)^{2\nu+1}}.
\end{equation*}
Thus, for $n\geq2$, we obtain
$$\max\{|a_n|,|b_n|\}\leq\frac{\beta^*_\nu(f)}{n}\sqrt{\frac{1+|\omega_f(0)|}{1-|\omega_f(0)|}}\frac{1}{r^{n-1}(1-r^2)^{\nu+1/2}}.$$
It is a simple exercise to see that $r^{1-n}(1-r^2)^{-(\nu+1/2)}$ is maximized in $r\in(0,1)$ for $r=\sqrt{\frac{n-1}{n+2\nu}}$.
Consequently,
\begin{align*}
\max\{|a_n|,|b_n|\}&\leq\frac{\beta^*_\nu(f)}{n}\sqrt{\frac{1+|\omega_f(0)|}{1-|\omega_f(0)|}}
\left(\frac{n+2\nu}{n-1}\right)^{n/2-1/2}\left(\frac{n+2\nu}{2\nu+1}\right)^{\nu+1/2}\\
&=\frac{\beta^*_\nu(f)\phi_\nu(n)}{(2\nu+1)^{\nu+1/2}}\sqrt{\frac{1+|\omega_f(0)|}{1-|\omega_f(0)|}}(n+2\nu)^{\nu-1/2},
\end{align*}
where $$ \phi_\nu(x)=\left[\left(1+\frac{2\nu+1}{x-1}\right)^{\frac{x-1}{2\nu+1}}\right]^{\nu+1/2}
\left(1+\frac{2\nu}{x}\right),\quad x\geq2.
$$

Next we prove that $\phi_\nu$ is an increasing function of $x$ to its limit $e^{\nu+1/2}$ in $[2,\infty)$.
Clearly, $\phi_\nu(x)>0$ for all $x\geq2$. For convenience, we let
$$\Phi_\nu(x)=(\log\phi_\nu(x))'=\frac{\phi'_\nu(x)}{\phi_\nu(x)}=\frac{1}{2}\log\left (\frac{x+2\nu}{x-1}\right )-\frac{(2\nu+1)x+4\nu}{2x(x+2\nu)}.
$$
Differentiating with respect to $x$ yields
$$\Phi'_\nu(x)
=-\frac{\psi_\nu(x)}{2x^2(x-1)(x+2\nu)^2}, \quad \psi_\nu(x)=(2\nu-1)^2x^2+8(\nu-\nu^2)x+8\nu^2 .
$$
If $\nu=1/2$, then $\psi_\nu(x)=2x+2\geq\psi_\nu(2)=6>0$ for all $x\geq2$.
If $\nu\neq1/2$, then we obtain
$$\psi'_\nu(x)=2(2\nu-1)^2x+8(\nu-\nu^2)\geq\psi'_\nu(2)=8(\nu-1/2)^2+2>0 ~\mbox{  for all $x\geq2$}
$$
and thus,  $\psi_\nu(x)\geq\psi_\nu(2)=4(2v^2+1)>0$ for all $x\geq2$.

Hence,  $\Phi'_\nu(x)<0$ in $[2,\infty)$ so that $\Phi_\nu(x)>\lim_{x\rightarrow\infty}\Phi_\nu(x)=0$ for all $x\in[2,\infty)$.
Therefore, we obtain $\phi'_\nu(x)>0$ in $[2,\infty)$ and the proof is complete.
\epf

\section{Bohr's inequalities} \label{sec5}

One of the classical problems in the theory of analytic functions which inspire many researchers is to determine
$$r_0=\sup \left \{r\in (0,1): \, M_f(r):=\sum_{n=0}^\infty |a_n|r^n\leq1\right \},
$$
where the supremum is taken over the class which consists of all functions of the form
$f(z)=\sum_{n=0}^\infty a_nz^n$ that converges in $\mathbb{D}$ and $|f(z)|\leq 1$ in $\mathbb{D}$.
It is well-known that $r_0=1/3$ and the number $1/3$ is called the classical Bohr radius for the class of all analytic self-maps
of the unit disk $\mathbb{D}$. Many authors have discussed the Bohr radius and extended this notion to various settings which led to
the introduction of Bohr's phenomenon.  As remarked in the introduction, we refer to \cite{AAP,KP,KPS} and  the references therein
for results on this topic.
 Moreover, in \cite{KP} the authors
introduced the notion of  $p$-Bohr radius for harmonic mappings which is defined as follows:
{\it  Let $f=h+\overline{g}$ be a harmonic mapping in $\mathbb{D}$, where $h$ and $g$ have the form \eqref{equ1a}.
For $p\geq1$, the $p$-Bohr radius for $f$ is defined to be the largest value $r_p$ such that
$$|a_0|+\sum_{n=1}^\infty (|a_n|^p+|b_n|^p)^{1/p}r^n\leq1 ~\mbox{ for }~ |z|=r\leq r_p.
$$
}
Clearly, all these radii coincide in the analytic case. The classical case $p=1$ is considered first time in \cite{abu}.

In this section, we determine the Bohr radius for analytic functions in $\mathcal{B}(\nu)$
and $p$-Bohr radius for harmonic mappings in $\mathcal{B}_H(\nu)$ and $\mathcal{B}^*_H(\nu)$.
The following results are  generalizations of that of the results of Kayumov  et al. \cite[Section~4]{KPS}.

\bthm \label{thm6}
Assume that $f(z)=\sum_{n=0}^\infty a_nz^n$ belongs to $\mathcal{B}(\nu)$ and $||f||_{\mathcal{B}(\nu)}\leq1$. Then
$$\sum_{n=0}^\infty |a_n|r^n\leq1
$$
for $|z|=r\leq r(\nu)=\max\{r_1(\nu),r_2(k)\}$
when $\nu\in(k/2,(k+1)/2]$ for some $k\in\mathbb{N}_0:=\mathbb{N}\cup \{0\}$.
Here $r_1(\nu)$ is the unique solution in $(0,1)$ to the equation
\beq \label{equ5}
6(1-r^2)^{2\nu}-\pi^2r^2=0
\eeq
and $r_2(k)$ is the unique solution in $(0,1)$  to the equation
\beq \label{equ6}
rF_k(r)-1+r=0,
\eeq
where
\begin{align} \label{equ7}
F_k(r)=
\begin{cases}
\ds \sum_{n=1}^{\infty}\frac{r^n}{n^2} &\mbox{for }~ k=0,\\
\ds \log\frac{1}{1-r} &\mbox{for }~ k=1,\\
\ds \frac{1}{k}\log\frac{1}{1-r}+\frac{1}{k}\sum_{n=1}^{k-1}\frac{1}{n}\left(\frac{1}{(1-r)^n}-1\right) &\mbox{for }~ k\geq2.
\end{cases}
\end{align}
Moreover,  $r(\nu)$ can not be replaced by $r_3(\nu)$ when $\nu\geq1$,
where
$$ r_3(1)=0.624162,~\mbox { and }~  r_3(\nu)=\min\left\{0.624162,~~\sqrt{1-1/\sqrt[\nu-1]{2\nu-1}}\right\}~\mbox{ for} ~\nu>1.
$$
\ethm
\bpf
By hypothesis, we have $||f||_{\mathcal{B}(\nu)}\leq1$ which gives
$$ \left |\sum_{n=1}^\infty na_nz^{n-1}\right |^2= |f'(z)|^2\leq\frac{(1-|a_0|)^2}{(1-|z|^2)^{2\nu}},\quad z\in\mathbb{D}.
$$
Integrating the inequality over the circle $|z|=r$ yields
\beq \label{equ8}
\sum_{n=1}^\infty n^2|a_n|^2r^{2(n-1)}\leq\frac{(1-|a_0|)^2}{(1-r^2)^{2\nu}}.
\eeq
By the classical Cauchy--Schwarz inequality, we obtain
\begin{align*}
|a_0|+\sum_{n=1}^\infty |a_n|r^{n}\leq&|a_0|+\sqrt{\sum_{n=1}^\infty n^2|a_n|^2r^{2n}}\sqrt{\sum_{n=1}^\infty \frac{1}{n^2}}\\
\leq&|a_0|+\frac{(1-|a_0|)r}{(1-r^2)^\nu}\sqrt{\frac{\pi^2}{6}}\\
\leq& 1 ~\mbox{ for }~ r\leq r_1(\nu),
\end{align*}
where $r_1(\nu)$ is the unique solution in (0,1) to the equation of \eqref{equ5}.
In fact, for each $\nu\in(0,\infty)$, the function $r/(1-r^2)^\nu$ increases from 0 to $\infty$ in  $[0,1)$.

On the other hand, if $\nu\in(\frac{k}{2},\frac{k+1}{2}]$ for some $k\in\mathbb{N}_0$,
then it follows from  \eqref{equ8} that
\beq \nonumber
\sum_{n=1}^\infty n^2|a_n|^2r^{n-1}\leq\frac{(1-|a_0|)^2}{(1-r)^{k+1}}.
\eeq
Integrating the above inequality twice (with respect to $r$) yields
$$\sum_{n=1}^\infty |a_n|^2r^{n}\leq(1-|a_0|)^2F_k(r),
$$
where $F_k(r)$ is defined by \eqref{equ7}. Applying the Cauchy--Schwarz inequality again, we have
\begin{align*}
|a_0|+\sum_{n=1}^\infty |a_n|r^{n}\leq&|a_0|+\sqrt{\sum_{n=1}^\infty |a_n|^2r^{n}}\sqrt{\sum_{n=1}^\infty r^n}\\
\leq&|a_0|+(1-|a_0|)\sqrt{\frac{F_k(r)r}{1-r}}\\
\leq& 1 ~\mbox{ for }~  r\leq r_2(k),
\end{align*}
where $r_2(k)$  is the unique solution to the equation of \eqref{equ6}.
Note that both $F_k(r)$ and $r/(1-r)$ are strictly increasing in $[0,1)$.
Combining the two estimates yields the desired conclusion.

For the upper bound of $r(\nu)$, since $\mathcal{B}(\nu)\supseteq\mathcal{B}(1)$ for any $\nu\geq1$,
it follows from \cite[Theorem~9]{KPS} that $r(\nu)$ can not be replaced by $0.624162$ for $\nu\geq1$.
In addition, let's consider the function
$$f_\nu(z)=\frac{(1-z^2)^{1-\nu}-1}{2(\nu-1)}=\sum_{n=1}^\infty a_{\nu,n}z^n,\quad z\in\mathbb{D}.
$$
A basic computation shows that $f_\nu\in\mathcal{B}(\nu)$ and $||f_\nu||_{\mathcal{B}(\nu)}=1$ when $\nu>1$.
It is easy to see that all coefficient $a_{\nu,n}$ are non-negative real number for each $\nu>1$ and $a_{\nu,n}=0$
for odd integer values of $n\geq 1$. For $\nu>1$, we consider the following inequality
$$
 \sum_{n=1}^\infty a_{\nu,n}r^n=\frac{(1-r^2)^{1-\nu}-1}{2(\nu-1)}\leq1,
$$
provided $r\leq \sqrt{1-1/\sqrt[\nu-1]{2\nu-1}}$ and thus, the conclusion follows.
\epf

It is easy to see that the function $r_1(\nu)$ is monotonically decreasing to 0 in $(0,+\infty)$.
In the following table, the notation $(r_1(\nu_1)\searrow r_1(\nu_2)]$ means that the value of $r_1(\nu)$
is monotonically decreasing from $\lim_{\nu\rightarrow\nu^+_1}r_1(\nu)=r_1(\nu_1)$ to $r_1(\nu_2)$ when $\nu_1<\nu\leq\nu_2$. 
So does $(r(\nu_1)\searrow r(\nu_2)]$.
Note that the function $r_3(\nu)$ is monotonically decreasing from $0.624162$ to 0 in $[1,+\infty)$ and $r_3(5.772240)=0.624162$.

\begin{table*}[htbp]


\begin{tabular}{cccc}

\toprule

$\nu$  & $r_1(\nu)$ & $r_2(k)$  & $r(\nu)$ \\

\midrule

(0,1/2] & $(0.779697\searrow0.614883]$ & $0.586028$ &  $(0.779697\searrow0.614883]$\\

(1/2,1] & $(0.614883\searrow0.546679]$ & $0.553567$ &  $(0.614883\searrow0.553567]$\\

(1,3/2] & $(0.546679\searrow0.503190]$ & $0.522089$ &  $(0.546679\searrow0.522089]$\\

(3/2,2] & $(0.503190\searrow0.471528]$ & $0.492552$ &  $(0.503190\searrow0.492552]$\\

(2,5/2] & $(0.471528\searrow0.446818]$ & $0.465403$ &  $(0.471528\searrow0.465403]$\\

(5/2,3] & $(0.446818\searrow0.426678]$ & $0.440723$ &  $(0.446818\searrow0.440723]$\\

\bottomrule

\end{tabular}

\end{table*}

\bthm \label{thm7}
Suppose that $f=h+\overline{g}\in\mathcal{B}_{H}(\nu)$,  where $h$ and $g$ are given by \eqref{equ1a}. If $||f||_{\mathcal{B}_{H}(\nu)}\leq1$
and $p\geq 1$, then we have
$$|a_0|+\sum_{n=1}^\infty (|a_n|^p+|b_n|^p)^{1/p}r^n\leq1
$$
for $|z|=r\leq\max\{r_{1}(\nu,p),r_{2}(k,p)\}$
when $\nu\in(k/2,(k+1)/2]$ for some $k\in\mathbb{N}_0$.
Here $r_{1}(\nu,p)$ is the unique solution in $(0,1)$ to the equation
\beq  \label{equ9}
6(1-r^2)^{2\nu}-M_p\pi^2r^2=0
\eeq
and  $r_{2}(k,p)$ is the unique solution in $(0,1)$ to the equation $M_p r F_k(r)-1+r=0,$
where $M_p=\max\{2^{(2/p)-1},1\}$ and $F_k(r)$ is defined by \eqref{equ7}.
\ethm

\bpf
By assumption, we see that
$$|h'(z)|^2+|g'(z)|^2\leq(|h'(z)|+|g'(z)|)^2\leq\frac{(1-|a_0|)^2}{(1-|z|^2)^{2\nu}},\quad z\in\mathbb{D}.
$$
Integrating the inequality over the circle $|z|=r$ so we get
\beq \nonumber
\sum_{n=1}^\infty n^2(|a_n|^2+|b_n|^2)r^{2(n-1)}\leq\frac{(1-|a_0|)^2}{(1-r^2)^{2\nu}}.
\eeq
Using the Cauchy-Schwarz inequality, we obtain
\begin{align*}
|a_0|+\sum_{n=1}^\infty (|a_n|^p+|b_n|^p)^{1/p}r^n
\leq&|a_0|+\sqrt{\sum_{n=1}^\infty n^2(|a_n|^p+|b_n|^p)^{2/p}r^{2n}}\sqrt{\sum_{n=1}^\infty \frac{1}{n^2}}\\
\leq&|a_0|+\sqrt{M_p\sum_{n=1}^\infty n^2(|a_n|^2+|b_n|^2)r^{2n}}\sqrt{ \frac{\pi^2}{6}}\\
\leq&|a_0|+\sqrt{M_p}\frac{(1-|a_0|)r}{(1-r^2)^\nu}\left (\frac{\pi}{\sqrt{6}}\right )
\end{align*}
which is less than or equal to $1$ provided
$r\leq r_{1}(\nu,p)$, where $r_{1}(\nu,p)$ is defined by \eqref{equ9}.
If $\nu\in(k/2,(k+1)/2]$ for some $k\in\mathbb{N}_0$, then we can combine the above proof with the corresponding proof of Theorem \ref{thm6}.
The resulting discussion completes the proof.
\epf

Next we will study $p$-Bohr radius for functions in $\mathcal{B}^*_H(\nu)$. Consider
$$f(z)=h(z)+\overline{g(z)}=\frac{1}{1-z}+\overline{\left (\frac{z}{1-z}\right )}
$$
so that $a_0=1$ and $a_n=b_n=1$ for $n\geq 1$.
Clearly, $f\in\mathcal{B}^*_H(\nu)$ and $||f||_{\mathcal{B}^*_H(\nu)}=|a_0|=1$ for any $\nu>0$.
However, we have $$|a_0|+\sum_{n=1}^\infty (|a_n|^p+|b_n|^p)^{1/p}r^n>|a_0|=||f||_{\mathcal{B}^*_H(\nu)}
$$ for any $r>0$.
In this case, the $p$-Bohr radius for $f$ is 0.
This is the reason why we add the condition of sense-preserving in the following result.

\bthm \label{thm8}
Suppose that $f=h+\overline{g}\in\mathcal{B}^*_H(\nu)$ is a sense-preserving harmonic mapping, where $h$ and $g$ are given by \eqref{equ1a}.
If $||f||_{\mathcal{B}^*_H(\nu)}\leq1$ and $p\geq1$, then
$$|a_0|+\sum_{n=1}^\infty (|a_n|^p+|b_n|^p)^{1/p}r^n\leq1
$$
for $|z|=r\leq\max\{r_{1}(\nu,p,|\omega_f(0)|),r_{2}(k,p,|\omega_f(0)|)\}$ when $\nu\in(k/2,(k+1)/2]$ for some $k\in\mathbb{N}_0$.
Here $r_{1}(\nu,p,|\omega_f(0)|)$ is the unique solution in $(0,1)$ to the equation
\beq \nonumber
3(1-|\omega_f(0)|)(1-r^2)^{2\nu+1}-M_p\pi^2(1+|\omega_f(0)|)r^2=0
\eeq
 and  $r_{2}(k,p,|\omega_f(0)|)$ is the unique solution in $(0,1)$ to the equation
\beq \nonumber
2M_p(1+|\omega_f(0)|)r F_{k+1}(r)-(1-|\omega_f(0)|)(1-r)=0,
\eeq
where $M_p=\max\{2^{2/p-1},1\}$ and $F_k(r)$ is defined by \eqref{equ7}.
\ethm
\bpf
By hypothesis $|g'(z)|< |h'(z)|$ and $||f||_{\mathcal{B}^*_H(\nu)}\leq1$ and thus, it follows from \eqref{equ4} that
$$|h'(z)|\leq \frac{1-|a_0|}{(1-|z|^2)^{\nu+\frac{1}{2}}}
\sqrt{\frac{1+|\omega_f(0)|}{1-|\omega_f(0)|}},\quad z\in\mathbb{D}.
$$
As in the proof of the previous theorem, we obtain that
$$\sum_{n=1}^\infty n^2(|a_n|^2+|b_n|^2)r^{2(n-1)}
\leq 2\frac{1+|\omega_f(0)|}{1-|\omega_f(0)|}\frac{(1-|a_0|)^2}{(1-r^2)^{2\nu+1}}.
$$
The remaining part of the proof is identical to Theorem \ref{thm7} and thus, we omit the details.
The proof is complete.
\epf

The dependence of $|\omega_f(0)|$ about $p$-Bohr radius in Theorem \ref{thm8}
can be seen from the following example.

\begin{example} \label{exa3}
{\rm For $ t\in[1/2,1)$, we consider the one parameter family of functions $f_t$ on $\ID$ given by
$$f_t(z)=h_t(z)+\overline{g_t(z)}=\sum_{n=0}^\infty a_{n}^{(t)} z^n+\overline{\sum_{n=1}^\infty b_{n}^{(t)} z^n},\quad z\in\mathbb{D},
$$
where
$$h_t(z)=1-2\sqrt{t-t^2}+\frac{1}{2}\log\frac{1+z}{1-z}\quad \text{and}\quad
g_t(z)=\frac{t-1}{2}\log(1-z^2)+\frac{t}{2}\log\frac{1+z}{1-z}.
$$
It is easy to see that each $f_t$ is sense-preserving in $\mathbb{D}$ with the dilatation $\omega_{f_t}(z)=(1-t)z+t$.
We find that
$$h_t'(z) =\frac{1}{1-z^2}, ~
|g_t'(z)| =\left|\frac{(1-t)z+t}{1-z^2}\right|\geq\frac{t-(1-t)|z|}{|1-z^2|}\geq\frac{2t-1}{|1-z^2|},\quad z\in\mathbb{D},
$$
and thus,
$$(1-|z|^2)\sqrt{J_{f_t}(z)}\leq(1-|z|^2)\sqrt{\frac{1}{|1-z^2|^2}-\frac{(2t-1)^2}{|1-z^2|^2}}\leq2\sqrt{t-t^2},
\quad z\in\mathbb{D},
$$
which implies that $f_t\in\mathcal{B}^*_H(1)$. Also, we observe that
$$(1-|x|^2)\sqrt{J_{f_t}(x)}
\rightarrow2\sqrt{t-t^2}~\mbox{ as }~ (-1,0)\ni x\rightarrow-1^+,
$$
which infers $\beta_1^*(f_t)=2\sqrt{t-t^2}$ and $||f_t||_{\mathcal{B}^*_H(1)}=1$.
Clearly, 
$$|a_0^{(t)}|+\sum_{n=1}^\infty (|a_n^{(t)}|^p+|b_n^{(t)}|^p)^{1/p}r^n>|a_0^{(t)}|=1-2\sqrt{t-t^2}
$$
for any $r>0$. Note that $|\omega_{f_t}(0)|=|g'_t(0)|/|h'_t(0)|=t$ and $1-2\sqrt{t-t^2}\rightarrow1=||f_t||_{\mathcal{B}^*_H(1)}$ as $t\rightarrow1^-$.
This means that if $t$ approaches to $1^{-}$, then the $p$-Bohr radius for  $f_t$ approaches  to $0$.
}
\end{example}

\subsection*{Acknowledgments}
The work was completed during the visit of the first author to the Indian Statistical Institute,
Chennai Centre and this author thanks the institute for the support and the hospitality.
The research of the first author was supported by the NSFs of China~(No. 11571049),
the Construct Program of the Key Discipline in Hunan Province and
the Science and Technology Plan Project of Hunan Province (No. 2016TP1020).
The second author is on leave from IIT Madras, Chennai.


\end{document}